\documentclass[11pt,openany]{article}

\usepackage [latin1]{inputenc}
\usepackage{amssymb}
\usepackage{amsmath}
\usepackage{amsthm}
\usepackage{amsfonts}
\usepackage{mathrsfs}
\usepackage{makeidx}
\usepackage{multicol}

\input{cyracc.def}

\numberwithin{equation}{section}

\begin{document}
%\pagenumbering{roman}
%\thispagestyle{empty}
\title{\Large \bf  Meromorphic functions partially share three values with their difference operators}

\author{Feng L\"{u}\footnote{Corresponding author, Email: lvfeng18@gmail.com, S19090058@s.upc.edu.cn.} ~ and Zhenliu Yang\\
%ENDAName
\small{College of Science, China University of Petroleum, Qingdao, Shandong, 266580, P.R. China.}\\
}
\date{}
\maketitle
%\vspace{1cm}

\vspace{3mm}

%\tableofcontents
%\newpage
%\centerline{\bf \large Abstract}
%\vspace{10mm}
\begin{abstract}
In this paper, we give a simple proof and strengthening of a uniqueness theorem of meromorphic functions which partially share 0, $\infty$ CM and 1 IM with their difference operators. Meanwhile, we partially solve a conjecture given by Chen-Yi in \cite{CY} and generalize some previous theorems in \cite{C, CX}.
\end{abstract}



\pagenumbering{arabic}

\section{Introduction and main results}

In this paper, the term "meromorphic function" always means meromorphic
in the whole complex plane $\mathbb{C}$. We say two meromorphic functions $f$ and $g$ share a constant $a$ IM if $f-a$ and $g-a$ have the same zeros or $ \overline{E}(a,f)=\overline E(a,g)$, and the notation $\overline{E}(a,f)$ denotes the set of all zeros of $f(z)-a$, where a zero is counted one time. In addition, we say $f$ and $g$ share $a$ CM if $f-a$ and $g-a$ have the same zeros with multiplicities or $E(a,f)= E(a,g)$, and the notation $E(a,f)$ denotes the set of all zeros of $f(z)-a$, where a zero with multiplicity $m$ is counted $m$ times. \\

The paper mainly concerns the uniqueness problem of meromorphic function $f(z)$ partially sharing 2 CM + 1 IM with its difference operator $\Delta_c f=f(z+c)-f(z)$, where $c$ is a nonzero constant. This problem was inspired by a conjecture given by Chen and Yi in \cite{CY}. They considered the uniqueness problem of meromorphic function $f(z)$ sharing $e_1$, $e_2$ and $\infty$ CM with $\Delta_c f$ and posed the following conjecture.\\

\noindent \textbf{Conjecture.} Let $f$ be a transcendental meromorphic function, and let $c\in \mathbb{C}\backslash \{0\}$ such that $\Delta_c f\not\equiv0$. If $\Delta_c f$ and $f(z)$ share three distinct values $e_1$, $e_2$, $\infty$ CM, then  $\triangle_cf=f$.\\

Since then, many scholars devoted to studying this conjecture. Zhang and Liao \cite{ZL} affirmed the conjecture for the case that $f$ is an entire function of finite order. Here and throughout this paper, the order and hyper order are defined in turn as follows:
$$
\rho (f)=\limsup_{r \rightarrow\infty} \frac{\log T(r,f)}{\log r},~~~~~\rho_{2}(f)=\limsup_{r \rightarrow\infty} \frac{\log\log T(r,f)}{\log r}.
$$

Later on, L\"{u} and L\"{u} \cite{LL} proved that the conjecture was also right when $f$ is a meromorphic function of finite order. Some related theorems can be found in \cite{Cui}. Recently, Chen and his co-workers considered the conjecture in another direction. They asked whether the shared value conditions ``3CM" in Conjecture can be weaken or not? Actually, they partially answered the question by the following two theorems, which were given by Chen in \cite{C} and  Chen and Xu in \cite{CX}, respectively.\\

\noindent \textbf{Theorem A.} Let $f$ be a nonconstant meromorphic function of hyper order $\rho_2(f)<1$, and let $c\in \mathbb{C}\backslash \{0\}$ such that $\Delta_c f\not\equiv 0$. If $\Delta_c f$ and $f(z)$ share 1 CM and satisfy $E(0,f(z))\subseteq E(0,\Delta_cf(z))$ and $E(\infty, \Delta_cf(z))\subseteq E(\infty, f(z))$, then  $\Delta_cf=f$.\\

\noindent \textbf{Theorem B.} Let $f$ be a nonconstant meromorphic function of hyper order $\rho_2(f)<1$, and let $c\in \mathbb{C}\backslash \{0\}$ such that $\Delta_c f\not\equiv 0$. If $\Delta_c f$ and $f(z)$ share 0, $\infty$ CM and 1 IM, then  $\Delta_c f=f$.\\

In this present paper, we still pay attention to this conjecture and give a simple proof and strengthening of a uniqueness theorem as follows.\\

\noindent \textbf{Theorem 1.} Let $f$ be a transcendental meromorphic function of hyper order $\rho_2(f)<1$, let $c\in \mathbb{C}\backslash \{0\}$, and let $L(f)=af(z+c)-bf(z)$, where $a,~b$ are two nonzero constants. If $E(0,f(z))\subseteq E(0, L(f))$, $E(\infty, L(f))\subseteq E(\infty, f(z))$ and $\overline{E}(1,f(z))\subseteq \overline{E}(1, L(f))$, then $L(f)=f(z)$ or $f(z+c)=f(z)$.\\

Obviously, Theorem 1 is an improvement of Theorems A and B. More generally, we deduce the following result and Theorem 1 is an immediate consequence of it.\\

\noindent \textbf{Theorem 2.} Let $f$ be a nonconstant meromorphic function of hyper order $\rho_2(f)<1$. If
$$
E(0, f(z))\subseteq E(0, f(z+c)),~~ E(\infty, f(z+c)) \subseteq E(\infty, f(z)),~~ \overline{E}(1, f(z)) \subseteq \overline{E}(A, f(z+c)),
$$
where $A$ is a nonzero constant, then $f(z+c)=Af(z)$ or $f(z+c)=f(z)$. In particularly, if $A=1$, then $f(z+c)=f(z)$.\\

\noindent \textbf{Remark 1.} We give the following examples to show that either $f(z+c)=Af(z)$ or $f(z+c)=\pm f(z)$ may occur. \\

 \textbf{Example 1.} Set $f(z)=e^{az}$, where $a$ is a nonzero constant such that $e^{ac}=A$. Then, $f(z)$ and $f(z+c)$ satisfies all the condition of Theorem 2. Clearly, $f(z+c)=Af(z)$. \\

 \textbf{Example 2.} Set $f(z)=1+e^{az}$, where $a$ is a nonzero constant such that $e^{ac}=1$. Then, $f(z+c)=f(z)$ and $1$ is a Picard value of $f(z)$. Obviously, $f(z)$ and $f(z+c)$ satisfies all the condition of Theorem 2.\\

\noindent \textbf{Remark 2.} Below, we apply Theorem 2 to derive Theorem 1. We write $f(z+c)$ as $f(z+c)=\frac{L(f)+bf(z)}{a}$. Then, the condition $E(0,f(z))\subseteq E(0,L(f)$ and $\overline{E}(1,f(z))\subseteq \overline{E}(1, L(f))$ yields that $E(0, f(z))\subseteq E(0, f(z+c))$ and $\overline{E}(1, f) \subseteq \overline{E}(\frac{1+b}{a}, f(z+c))$, respectively. Set $f(z)=\frac{af(z+c)-L(f)}{b}$. Suppose that $z_0$ is a pole of $f(z+c)$ with multiplicity $m$. We assume that $z_0$ is a pole of $f(z)$ with multiplicity $n$ (It is pointed out that $n$ maybe zero). Below, we prove $m\leq n$. On the contrary, assume $n<m$. Then $z_0$ is a pole of $L(f)=af(z+c)-bf(z)$ with multiplicity $m$. The condition $E(\infty, L(f))\subseteq E(\infty, f(z))$ yields that $z_0$ is a pole of $f(z)$ with multiplicity at least $m$, a contradiction. Thus $m\leq n$ and $E(\infty, f(z+c)) \subseteq E(\infty, f(z))$. The above discussions yields that $f(z)$ and $f(z+c)$ satisfies all the conditions of Theorem 2. So, $f(z+c)=\frac{1+b}{a}f(z)$ or $f(z+c)=f(z)$. The first case implies that $L(f)=f$. Therefore, the proof of Theorem 1 is finished.\\

Before to proceed, we spare the reader for a moment and assume his/her familiarity with
the basics of Nevanlinna's theory of meromorphic functions in $\mathbb{C}$ and the usual notations such as
$T(r,f)$, $m(r,f)$, $N(r,f)$, $\overline{N}(r,f)$ and $S(r,f)$ (see e.g., \cite{WK, CCY}). \\

\begin{proof}[Proof of Theorem 2] The proof is based on some idea of Chen in \cite{C1}. On the contrary, we assume that $f(z+c)\not\equiv Af(z)$ and $f(z+c)\not\equiv f(z)$. Below, we will derive a contradiction. We firstly introduce the auxiliary function
\begin{equation}\label{2.0}
H(z)=\frac{f(z+c)}{f(z)}.
\end{equation}
By $\rho_2(f)<1$ and the difference analogue of the logarithmic derivative lemma (which can be found in \cite{HAL}), one can easily get $m(r,H)=S(r,f)$. The conditions $E(0, f(z)) \subset E(0, f(z+c))$ and $E(\infty, f(z+c)) \subset E(\infty, f(z))$ yields that
$H$ is an entire function. Thus, $T(r,H)=S(r,f)$. We claim that $\overline{N}(r,\frac{1}{f(z)-1})=S(r,f)$.\\

If 1 is a Picard value of $f(z)$, then the claim is right. Below, assume that 1 is not a Picard value of $f(z)$ and $z_{0}$ is a zero of $f(z)-1$. Then, $f(z_0+c)=A$ and $H(z_0)=A$. If $H(z)\equiv A$, then $f(z+c)=Af(z)$, a contradiction. So, $H(z)\not\equiv A$. Thus,
\begin{equation}\label{2.1}
\overline{N}(r,\frac{1}{f(z)-1})\leq  \overline{N}(r,\frac{1}{H(z)-A})\leq T(r,H)+O(1)=S(r,f),
\end{equation}
which implies that the claim holds. Below, we employ a result, which is Lemma 8.3 of Halburd-Korhonen-Tohge in \cite{HAL}.\\

\emph{\noindent{\bf Lemma.} Let $T:[0,+\infty) \rightarrow[0,+\infty)$ be a non-decreasing continuous function and let $s \in(0, \infty) .$ If the hyper-order of $T$ is strictly less than one, i.e.,
$$\limsup _{r \rightarrow \infty} \frac{\log \log T(r)}{\log r}=s<1$$
and $\delta \in(0,1-\mathrm{s}),$ then
$$T(r+s)=T(r)+o\left(\frac{T(r)}{r^{\delta}}\right)$$
where $r$ runs to infinity outside of a set of finite logarithmic measure.}\\

Applying the above lemma to $f(z)$, one gets
\begin{equation}\label{2.2}
\overline{N}(r,\frac{1}{f(z)-1})=\overline{N}(r,\frac{1}{f(z-c)-1})+S(r,f)=\overline{N}(r,\frac{1}{f(z+c)-1})+S(r,f).
\end{equation}¡¡¡¡
Therefore,
\begin{equation}\label{2.3}
\overline{N}(r,\frac{1}{f(z-c)-1})=S(r,f),~~\overline{N}(r,\frac{1}{f(z+c)-1})=S(r,f).
\end{equation}
We rewrite (\ref{2.0}) as
\begin{equation}\label{2.4}
f(z+c)-1=H(z)[f(z)-\frac{1}{H(z)}].
\end{equation}
Then, combining (\ref{2.3}) and (\ref{2.4}) yields that
\begin{equation}\label{2.5}
\overline{N}(r,\frac{1}{f(z)-\frac{1}{H(z)}})=\overline{N}(r,\frac{1}{f(z+c)-1})+S(r,f)=S(r,f).
\end{equation}
We rewrite (\ref{2.0}) as $f(z)=H(z-c)f(z-c)$. Then¡¡¡¡¡¡¡¡
\begin{equation}
f(z-c)-1=\frac{1}{H(z-c)}[f(z)-H(z-c)].
\end{equation}
The same argument leads to
\begin{equation}\label{2.6}
\overline{N}(r,\frac{1}{f(z)-H(z-c)})=\overline{N}(r,\frac{1}{f(z-c)-1})+S(r,f)=S(r,f).
\end{equation}
Assume that the functions $\frac{1}{H(z)}$, $H(z-c)$ and 1 are distinct with each other. Applying the second main theorem of Nevanlinna, one gets
$$
T(r,f(z))\leq \overline{N}(r,\frac{1}{f(z)-H(z-c)})+\overline{N}(r,\frac{1}{f(z)-\frac{1}{H(z)}})+\overline{N}(r,\frac{1}{f(z)-1})+S(r,f)=S(r,f),
$$
a contradiction. Therefore, there are at least two functions in the set $\{\frac{1}{H(z)}, ~H(z-c),~~1\}$ which are equal identically. If $\frac{1}{H(z)}=1$ or $H(z-c)=1$, then $f(z)\equiv f(z+c)$, a contradiction. Therefore, it is suffice to handle the case $\frac{1}{H(z)}=H(z-c)$. \\

We rewrite the above equation as $H(z)H(z+c)=1$. Then, $H^2(z)=\frac{H(z)}{H(z+c)}$. Suppose that $H$ is not a constant. Then,
$$
2T(r,H(z))=2m(r,H(z))=m(r,H^2(z))=m(r,\frac{H(z)}{H(z+c)})=S(r,H),
$$
a contradiction. Thus, $H$ is a constant. The above equation yields that $H=\pm 1$. So, $f(z+c)=-f(z)$. If $A=-1$, then $f(z+c)=Af(z)$, a contradiction, which implies that $A\neq -1$. Further, the condition $\overline{E}(1, f(z)) \subset \overline{E}(A, f(z+c))$ yields that 1 is a Picard value of $f(z)$. Meanwhile, $-1$ is a Picard value of $f(z+c)$, so is $f(z)$. Thus, we can set $\frac{f(z)-1}{f(z)+1}=e^\beta$, where $\beta$ is a nonconstant entire function with order $\rho(\beta)<1$. We rewrite the above equation as $f(z)=\frac{1+e^\beta}{1-e^\beta}$. Substituting the form of $f(z)$ into the function $f(z+c)=-f(z)$ yields that $e^{\beta(z+c)+\beta(z)}=1$, which leads to $\beta(z+c)+\beta(z)=2k\pi i$, $k$ is a fixed integer. Furthermore, one gets $\rho(\beta)\geq 1$, a contradiction. \\

Therefore, we finish the proof of this theorem.\\

\end{proof}

To conclude this paper, we give two natural further studies which
are related with the main results. One is to generalize the function $L(f)$ in Theorem 1 to higher order linear difference polynomial $L_n(f)=\sum_{j=1}^n a_jf(z+jc)$, where $a_j$'s are constant, or more generally small functions with respect to $f(z)$. For special case, $L_n(f)$ becomes to the $n$-th difference operator $\Delta_c ^nf$. The other one is to weaken the condition $E(a,f(z))\subseteq E(a, L(f))$ in Theorem 1 to $E(a,f(z))\setminus G\subseteq E(a, L(f))$, where $a$ is a constant and $G$ is an set. And the set $G$ is called an exceptional set. One would like the exceptional set $G$ as large as possible.

\end{document}